\documentclass{amsart}

\usepackage{graphicx}

\usepackage{pst-node}
\usepackage{tikz-cd}

\theoremstyle{definition}

\theoremstyle{remark}

\numberwithin{equation}{section}

\begin{document}

\title{Alexander and Jones polynomials of surgerized tst links}

\author{Wilson Wong, Franky Mok}
\address{Email: wilsonwys@yahoo.com.hk}


\date{January 31, 2019}

\dedicatory{Dedicated to our team members Sam Lau, Aaron Chan and Kenith Chan, and our teacher Mr. S. H. Chan, with whom we worked together on a HLMA project in 2012.}


\begin{abstract}
This paper is a continuation on the 2012 paper on "Cutting Twisted Solid Tori (TSTs)", in which we considered twisted solid torus links (tst links). We generalize the notion of tst links to "surgerized tst links": recall that when performing $\Phi^\mu(n(\tau), d(\tau), M)$ on a tst $\langle \tau \rangle$ where $M$ is odd, we obtain the tst link, $[\Phi^\mu(n(\tau), d(\tau), M)]$ that contains a trivial knot as one of its components. We then perform another operation  $\Phi^{\mu'}(n(\tau  '), d(\tau '), M')$ on that trivial knot to create a new link, which we call a "surgerized tst link" (stst link). If $M'$ is odd, we can repeat the process to give more complicated stst links. We compute braid words, Alexander and Jones polynomials of such links.
\end{abstract}

\maketitle
\section*{Introduction}
Recall that a twisted solid torus link (tst link) $[\Phi^\mu(n(\tau), d(\tau), M)]$ is obtained by cutting a twisted solid torus (tst) $\langle \tau \rangle$ into $M$ pieces with chosen denominator $d(\tau)$ by inserting $\mu(M-1)$ blades. For details please see \cite{HLMA}. The resultant link has $(\mu P + I)$ components, where $P$ and $I$ are defined in the paper \cite{HLMA} using the number of partial tsts attached to different faces in the fundamental graph $\Lambda(n(\tau), d(\tau), M)$. When $N(\tau)\neq \tau$, $P$ and $I$ are respectively the same as the multiplicity of knotted tst $\langle N(\tau) \rangle_K$ and $\langle \tau\rangle_K$ in the knotted tst sum $\Phi(n(\tau), d(\tau), M)$. When $N(\tau) = \tau$, the multiplicity of $\langle \tau\rangle_K$ is $P + I$. In particular $I$ has a simple formula: 
$$I=\dfrac{1-(-1)^M}2 = M\mod 2$$
This means that $I$ is $1$ when $M$ is odd and $0$ when $M$ is even. The tst link  $[\Phi^\mu(n(\tau), d(\tau), M)]$ is a $(-\mu P N(\tau), \mu P D(\tau))$-torus link when $M$ is even, where each component is a $\langle N(\tau) \rangle_K$, and a $(-\mu P N(\tau), \mu P D(\tau))$-torus link with a trivial knot, the components being $\langle N(\tau) \rangle_K$ and $\langle \tau\rangle_K$ respectively when $M$ is odd.

Thus, if $M$ is odd, we can consider the following surgery: first, we drill $\langle \tau\rangle_K$ from $S^3$ with the $(-\mu P N(\tau), \mu P D(\tau))$-torus link left untouched. Then we replace the part that is drilled away with a new twisted solid torus $\langle \tau'\rangle$. Afterwards we perform the operation $\Phi^{\mu'}(n(\tau'), d(\tau'), M')$ on this new tst. 

The resultant knotted tst sum is then denoted $\Phi(n(\tau),d(\tau), M, \mu; n(\tau'),d(\tau'),M',\mu')$, and the resultant link is denoted $[\Phi(n(\tau),d(\tau), M, \mu; n(\tau'),d(\tau'),M',\mu')]$ and called a surgerized tst link (stst link). If we have
$$\Phi^{\mu'}(n(\tau'), d(\tau'), M') = \mu'P' \langle N(\tau')\rangle_K + I' \langle \tau'\rangle_K $$
then $[\Phi(n(\tau),d(\tau), M, \mu; n(\tau'),d(\tau'),M',\mu')]$ has $(\mu P + \mu'P'+ I')$ components, and
$$\Phi(n(\tau),d(\tau), M, \mu; n(\tau'),d(\tau'),M',\mu') =\mu P \langle N(\tau)\rangle_K  + \mu'P' \langle N(\tau')\rangle_K + I' \langle \tau'\rangle_K $$
If $M'$ is odd, we can repeat the process and obtain the stst link $$[\Phi(n(\tau),d(\tau), M, \mu;n(\tau'),d(\tau'), M', \mu'; n(\tau''),d(\tau''), M'', \mu'')]$$ and so forth. Performing the surgery $(k-1)$ times and performing the cut $\Phi^{ \mu_j}(n(\tau_j),d(\tau_j), M_j)$ the $(j-1)^\text{st}$ time, where $M_1,M_2, ..., M_{k-1}$ are odd, we obtain the stst link $$[\Phi(n(\tau_1),d(\tau_1), M_1, \mu_1;n(\tau_2),d(\tau_2), M_2, \mu_2; ...; n(\tau_k),d(\tau_k), M_k, \mu_k)]$$
or more compactly denoted by $[\Phi_{j=1}^k(n(\tau_j),d(\tau_j), M_j, \mu_j)]$. We also have the knotted tst sum given by
$$\Phi_{j=1}^k(n(\tau_j),d(\tau_j), M_j, \mu_j) = \sum_{j=1}^k \mu_j P_j\langle N(\tau_j)\rangle_K + I_k\langle \tau_k\rangle_K$$
where  $\Phi^{\mu_j}(n(\tau_j),d(\tau_j), M_j) = \mu_j P_j \langle N(\tau_j)\rangle_K + \langle \tau_j\rangle_K$ for $j=1,2,...,k-1$, and $\Phi^{\mu_k}(n(\tau_k),d(\tau_k), M_k) = \mu_k P_k \langle N(\tau_k)\rangle_K + I_k\langle \tau_k\rangle_K$. The link $[\Phi_{j=1}^k(n(\tau_j),d(\tau_j), M_j, \mu_j)]$ has  $(\sum_{j=1}^k \mu_j P_j + I_k)$ components, and consists of $\mu_j P_j$ torus knots of the $(-N(\tau_j),D(\tau_j) )$ type, and $I_k$ trivial knot(s). 

In this paper, we will compute a braid word for each of the stst links, then a Seifert matrix using the braid words. We will then compute the Alexander polynomials using techniques of Heegaard splitting and Fox calculus
. Then, we will compute the Jones polynomials using techniques of topological quantum field theory and Heegaard splitting in a different manner.

\section{Braid words} 
Recall that a braid word of $[\Phi^\mu(n(\tau),d(\tau),M)]$ is $(\sigma_1^{-1}\sigma_2^{-1}\cdots \sigma_{\mu PD(\tau)-1}^{-1})^{\mu PN(\tau)}$ when $M$ is even, and $(\sigma_1^{-2}\sigma_2^{-1}\cdots \sigma_{\mu PD(\tau)}^{-1})^{\mu PN(\tau)}$ when $M$ is odd as we derived  in \cite{HLMA}. We now write $N(\tau_j) = N_j$,$D(\tau_j) =D_j$ and $\mu_jP_j=S_j$ and consider a braid word for the stst link $[\Phi_{j=1}^k(n(\tau_j),d(\tau_j), M_j, \mu_j)]$.

For $k=2$ and $M_2$ even, the braid word is given by
$$\left(\coprod_1 ^{i=S_2D_2} \sigma_i^{-1} \prod_{i=1} ^{S_1D_1+S_2D_2-1}\sigma_i^{-1} \right) ^{S_1N_1} \left (\prod_{i=1}^{S_2D_2-1} \sigma_i^{-1} \right) ^{S_2N_2} \in B_{S_1D_1+S_2D_2} $$
where $\prod_{i=a}^b \beta_i = \beta_a\beta_{a+1}\cdots \beta_b$, $\coprod_a^{i=b} \beta_i = \beta_b \beta_{b-1}\cdots \beta_a$, $\beta_i$ is a braid-valued function on $i$ and $B_n$ is the braid group of $n$ strings.

For $M_2$ odd, the braid word is given by
$$\left(\coprod_1 ^{i=S_2D_2+1} \sigma_i^{-1} \prod_{i=1} ^{S_1D_1+S_2D_2}\sigma_i^{-1} \right) ^{S_1N_1} \left (\sigma_1^{-1} \prod_{i=1}^{S_2D_2} \sigma_i^{-1} \right) ^{S_2N_2}\in B_{S_1D_1+S_2D_2+1}  $$
Generalizing, a braid word of the stst link for even $M_k$ is given by
\begin{equation*}
\begin{split}
&\left(\coprod_1 ^{i=\sum_{j=2}^k S_jD_j} \sigma_i^{-1} \prod_{i=1} ^{\sum_{j=1}^k S_jD_j-1}\sigma_i^{-1} \right) ^{S_1N_1} \left (\coprod_1 ^{i=\sum_{j=3}^k S_jD_j} \sigma_i^{-1} \prod_{i=1}^{\sum_{j=2}^k S_jD_j-1} \sigma_i^{-1} \right) ^{S_2N_2}\cdots \left (\prod_{i=1}^{S_kD_k-1} \sigma_i^{-1} \right) ^{S_kN_k}\\
=&\prod_{l=1}^k\left(\coprod_1 ^{i=\sum_{j=l+1}^k S_jD_j} \sigma_i^{-1} \prod_{i=1} ^{\sum_{j=1}^k S_jD_j-1}\sigma_i^{-1} \right) ^{S_lN_l} \in B_{\sum_{j=1}^k S_jD_j}
\end{split}
\end{equation*}
where $\sum_{j=k+1}^k$ is interpreted as $0$ and $\coprod_1^{i=0}$ is interpreted as the identity braid, whereas for odd $M_k$, it is given by
$$\prod_{l=1}^k\left(\coprod_1 ^{i=\sum_{j=l+1}^k S_jD_j+1} \sigma_i^{-1} \prod_{i=1} ^{\sum_{j=1}^k S_jD_j}\sigma_i^{-1} \right) ^{S_lN_l} \in B_{\sum_{j=1}^k S_jD_j+1}$$

\section{Alexander polynomials}
To find the Alexander polynomials of the stst link $[\Phi_{j=1}^k(n(\tau_j),d(\tau_j), M_j, \mu_j)]$, we first review how the Alexander polynomials of the tst link $[\Phi^\mu(n(\tau),d(\tau), M)]$ is computed.

\subsection{Tst links $[\Phi^\mu(n(\tau), d(\tau), M)]$ with even $M$}
For even integers $M$, the link is the same as the torus link $K_{-\mu P N(\tau), \mu P D(\tau)}$. Following \cite{Murasugi}, the Alexander polynomial of the link is given by
$$\Delta_{ [\Phi^\mu(n(\tau),d(\tau), M)]} (t) = (-1)^{\mu P -1}t^{-\frac{(\mu PN(\tau)-1)(\mu PD(\tau)-1)}2}\dfrac{(1-t)(1-t^{\mu P N(\tau)D(\tau)})^{\mu P}}{(1-t^{\mu P N(\tau)})(1-t^{\mu P D(\tau)})}$$
Murasugi claimed that this can be computed from the Seifert matrix of the link \cite{Murasugi}, but it is in fact very difficult to do so, and which is probably why the book skipped the derivation.

We now derive this result using Fox calculus and Seifert-van Kampen Theorem. For the sake of simplicity we write $\mu P = S$, $N(\tau) = N$ and $D(\tau) = D$ and consider the Alexander polynomial of the link $L =K_{SN, SD}$. Let $n(C)$ denote the tubular neighborhood of a link $C$ in $S^3$. Firstly we note that if we put all the $S$ knots onto the same torus $\partial U$ that bounds a solid torus $U$ and separate the complement of $L$ along the torus into the manifolds $M = \overline{S^3 \backslash (U\cup n(L))} $  and $N = \overline {U\backslash n(L)}$, their intersection $M\cap N$ is not connected, therefore Seifert-van Kampen Theorem cannot be applied (at least not in its common form). Instead, we consider $S$ solid tori $U_S \subset \cdots \subset U_2 \subset U_1$ that share the same longitudinal circle. On each of their boundaries $\partial U_i$, we place one torus knot $K_{N,D}$ and let $K_i$ denote that knot. Consider the $3$-manifolds $X_i = \overline {U_i \backslash (U_i\cap n(K_i))\backslash (U_{i+1}\cup n(K_{i+1})) }$ for $1\leq i \leq S-1$, $X_0 = \overline{S^3\backslash (U_1\cup n(K_1))}$ and $X_S = \overline{U_S\backslash (U_S\cap K_S)}$. Then $\bigcup_{i=0}^SX_i=\overline{S^3\backslash n(L)} = S^3\backslash L$. Now $X_{i-1} \cap X_i = \overline{\partial U_i \backslash (\partial U_i\cap n(K_i))}$ is connected, so we can apply Seifert-van Kampen Theorem. Let $a_i$ and $b_i$ be the meridian and longitudinal circles respectively inside the solid torus $U_i$ for $1\leq i \leq S$. Note that $a_S$ is contractible inside $X_S$. We also introduce its meridian circle of the solid torus $V = S^3\backslash U_1$ to be $a_0$. The longitudinal circle of $V$ is contractible in $X_0$. The fundamental group $\pi_1(X_0)$ of $X_0$ is generated by $a_0$, and that of $X_i$ is
$$\pi_1(X_i) \cong \langle a_i, b_i | a_ib_i = b_i a_i\rangle$$
where  $1\leq i \leq S-1$. Now we consider the fundamental group of $X_{i-1} \cap X_i$:
it is generated by $g_i$, which is the knot $K_i$ slightly translated, so that it sits on the intersection. Now note that $g_i$ winds inside $U_{i-1}$ around the longitudinal circle $D$ times, and $N$ times around the meridian circle, so we have $g_i=a_{i-1}^Nb_{i-1}^D$. Similarly we have $g = a_i^N b_i^D$ inside $U_i$. Repeatedly applying Seifert-van Kampen Theorem gives the fundamental group of the link complement:
$$\pi_1(S^3\backslash L)\cong \left \langle \begin{split} a_0, a_1, a_2, ..., a_{S - 1},\\ b_1, b_2, ..., b_S \end{split}\Biggr|\begin{split}
a_0^N = a_1^Nb_1^D, & a_i^Nb_i^D = a_{i+1}^Nb_{i+1}^D= b_S^D, &  1\leq i \leq S-2 \\
a_jb_j = b_j a_j,              & 1\leq j \leq S-1                        & \\
\end{split}\right\rangle$$
Applying a Tietze type II transformation gives
$$\pi_1(S^3\backslash L)\cong \left \langle \begin{split} a_0, a_1, a_2, ..., a_{S - 1},\\ b_1, b_2, ..., b_S \end{split}\Biggr|\begin{split}
a_0^N &= a_i^Nb_i^D =b_S^D,&  1\leq i \leq S-1\\
a_jb_j &= b_j a_j, & 1\leq j \leq S-1\\
\end{split}\right\rangle$$
Next we use Fox calculus to find the Alexander matrix. Recall that under the Fox differential operator we have
$$\dfrac{\partial (fg)}{\partial x} = \dfrac{\partial f}{\partial x} + f\dfrac{\partial g}{\partial x} $$
Thus we have
\begin{equation*}
\begin{split}
\dfrac{\partial x^n}{\partial x} &= \dfrac{\partial x^{n-1}}{\partial x} + x^{n-1}\dfrac{\partial x}{\partial x} \\
                                                        &= \dfrac{\partial x^{n-1}}{\partial x} + x^{n-1} \\
                                                        &= \dfrac{\partial x^{n-2}}{\partial x} + x^{n-2} + x^{n-1}\\
                                                        &\vdots \\
                                                        &=\dfrac{\partial x}{\partial x}+ x + \cdots +x^{n-2} + x^{n-1}\\
                                                        &=1+x+ \cdots + x^{n-1} \\
                                                        &=\dfrac{1-x^n}{1-x}
\end{split}
\end{equation*}
We then find the Fox derivative of an inverse note that $ff^{-1} = 1$ and $\frac{\partial (1) }{\partial x} = 0$:
\begin{equation*}
\begin{split}
\dfrac{\partial (f f^{-1})}{\partial x} &= \dfrac{\partial f}{\partial x} + f\dfrac{\partial f^{-1}}{\partial x} \\
0                                                              &= \dfrac{\partial f}{\partial x} + f\dfrac{\partial f^{-1}}{\partial x} \\
0                                                              &= f^{-1} \dfrac{\partial f}{\partial x} + \dfrac{\partial f^{-1}}{\partial x} \\
\dfrac{\partial f^{-1}}{\partial x}     &= -f^{-1} \dfrac{\partial f}{\partial x}
\end{split}
\end{equation*}
For the relation $f=g$, we can set up a relator $r = fg^{-1}$, for which the Fox derivative is given by
\begin{equation*}
\begin{split}
\dfrac{\partial r}{\partial x} &= \dfrac{\partial (fg^{-1})}{\partial x} \\
                                                    &= \dfrac{\partial f}{\partial x} + f\dfrac{\partial g^{-1}}{\partial x} \\
                                                    &= \dfrac{\partial f}{\partial x} - fg^{-1}\dfrac{\partial g}{\partial x} \\
                                                   &=\dfrac{\partial f}{\partial x} - (1)\dfrac{\partial g}{\partial x} \\
                                                   &= \dfrac{\partial f}{\partial x} -\dfrac{\partial g}{\partial x}
\end{split}
\end{equation*}
where in the penultimate step we have set $f=g$ so that $fg^{-1} = 1$.

Using these formulae we can then find the Fox derivatives for our relators. Let $r_i = a_0^N (a_i^N b_i^D)^{-1}$, $r_S = a_0^N b_S^{-D}$ and $s_j = (a_jb_j)(b_ja_j)^{-1}$. Then we have
\begin{equation*}
\begin{split}
\dfrac{\partial r_i}{\partial a_0} &= \dfrac{1-a_0^N}{1-a_0} = \dfrac{\partial r_S}{\partial a_0} \\
\dfrac{\partial r_i}{\partial a_i} &= -\dfrac{1-a_i^N}{1-a_i}\\
\dfrac{\partial r_i}{\partial b_i} &= -a_i^N\dfrac{1-b_i^D}{1-b_i} \\
\dfrac{\partial r_S}{\partial b_S} &= -\dfrac{1-b_S^D}{1-b_S} \\
\dfrac{\partial s_j}{\partial a_j} &= 1-b_j\\
\dfrac{\partial s_j}{\partial b_j} &= a_j -1 
\end{split}
\end{equation*}
where $1\leq i,j \leq S -1$. Other Fox derivatives are zero.
An Alexander matrix for our link $L$ is given by
$$\begin{bmatrix}
\frac{1-a_0^N}{1-a_0} & -\frac{1-a_1^N}{1-a_1} &&&& -a_1^N\frac{1-b_1^D}{1-b_1} \\
\frac{1-a_0^N}{1-a_0} & & -\frac{1-a_2^N}{1-a_2} &&&&-a_2^N\frac{1-b_2^D}{1-b_2} \\
\vdots & &&\ddots  &&&&\ddots\\
\frac{1-a_0^N}{1-a_0} & && &-\frac{1-a_{S -1}^N}{1-a_{S -1}} &&&&-a_{S -1}^N\frac{1-b_{S -1}^D}{1-b_{S -1}}\\
\frac{1-a_0^N}{1-a_0} & && &&&&&&-\frac{1-b_S^D}{1-b_S}\\
& 1-b_1&&&&a_1-1\\
&&1-b_2&&&&a_2-1 \\
&&&\ddots &&&&\ddots \\
&&&&1-b_{S -1}&&&&a_{S -1} -1
\end{bmatrix}$$
where unspecified entries are zero. The matrix is of size $(2S-1)\times 2S$.

The Alexander ideal is generated by determinants of the following types:

$\textbf{Type I:}$ Here the first column is deleted.
\begin{equation*}
\begin{split}
&\begin{vmatrix}
 -\frac{1-a_1^N}{1-a_1} &&&& -a_1^N\frac{1-b_1^D}{1-b_1} \\
 & -\frac{1-a_2^N}{1-a_2} &&&&-a_2^N\frac{1-b_2^D}{1-b_2} \\
 &&\ddots  &&&&\ddots\\
 && &-\frac{1-a_{S -1}^N}{1-a_{S -1}} &&&&-a_{S -1}^N\frac{1-b_{S -1}^D}{1-b_{S -1}}\\
&& &&&&&&-\frac{1-b_S^D}{1-b_S}\\
 1-b_1&&&&a_1-1\\
&1-b_2&&&&a_2-1 \\
&&\ddots &&&&\ddots \\
&&&1-b_{S -1}&&&&a_{S -1} -1
\end{vmatrix}\\
=&-(-1)^{2(S-1)+S}\dfrac{1-b_S^D}{1-b_S}\begin{vmatrix}
 -\frac{1-a_1^N}{1-a_1} &&&& -a_1^N\frac{1-b_1^D}{1-b_1} \\
 & -\frac{1-a_2^N}{1-a_2} &&&&-a_2^N\frac{1-b_2^D}{1-b_2} \\
 &&\ddots  &&&&\ddots\\
 && &-\frac{1-a_{S -1}^N}{1-a_{S -1}} &&&&-a_{S -1}^N\frac{1-b_{S -1}^D}{1-b_{S -1}}\\
 1-b_1&&&&a_1-1\\
&1-b_2&&&&a_2-1 \\
&&\ddots &&&&\ddots \\
&&&1-b_{S -1}&&&&a_{S -1} -1
\end{vmatrix}
\end{split}
\end{equation*}
\begin{equation*}
\begin{split}
=&(-1)^{S+1}\dfrac{1-b_S^D}{1-b_S}\left|
\begin{split}
\begin{bmatrix}
-\frac{1-a_1^N}{1-a_1}  \\
 & -\frac{1-a_2^N}{1-a_2} \\
 &&\ddots \\
 && &-\frac{1-a_{S -1}^N}{1-a_{S -1}} 
\end{bmatrix}
\begin{bmatrix}
a_1-1\\
&a_2-1 \\
&&\ddots \\
&&&a_{S -1} -1
\end{bmatrix}\\
-\begin{bmatrix}
 -a_1^N\frac{1-b_1^D}{1-b_1} \\
 &-a_2^N\frac{1-b_2^D}{1-b_2} \\
 &&\ddots\\
 && &-a_{S -1}^N\frac{1-b_{S -1}^D}{1-b_{S -1}}
\end{bmatrix}
\begin{bmatrix}
 1-b_1\\
&1-b_2 \\
&&\ddots  \\
&&& 1-b_{S -1}
\end{bmatrix}
\end{split}
\right|\\
=&(-1)^{S+1}\dfrac{1-b_S^D}{1-b_S} 
\begin{vmatrix}
1-a_1^N+ a_1^N - a_1^Nb_1^D  \\
 & 1-a_2^N+ a_2^N - a_2^Nb_2^D \\
 &&\ddots \\
 && &1-a_{S-1}^N+ a_{S-1}^N - a_{S-1}^Nb_{S-1}^D 
\end{vmatrix}\\
=&(-1)^{S+1}\dfrac{(1-b_S^D)\prod_{i=1}^{S-1}(1-a_i^Nb_i^D)}{1-b_S}\\
=&(-1)^{S+1}\dfrac{(1-a_0^N)^S}{1-b_S}
\end{split}
\end{equation*}
where we have used Schur complement formula and that the bottom two matrix blocks are diagonal matrices and therefore commute, and set the relators equal to one.
\clearpage
$\textbf{Type II.j:}$ Here the $(j+1)^\text{st}$ column is deleted, where $1\leq j \leq S-1$. We exchange columns so that the $j^\text{th}$ column in the Type I determinant is replaced by the first column of the Alexander matrix.

\begin{equation*}
\begin{split}
&(-1)^j\begin{vmatrix}
 -\frac{1-a_1^N}{1-a_1} && &\frac{1-a_0^N}{1-a_0}  \\
 & -\frac{1-a_2^N}{1-a_2} &&\frac{1-a_0^N}{1-a_0} \\
 &&\ddots  &\vdots \\
 &&&\frac{1-a_0^N}{1-a_0}&&&\textbf{A} \\
 &&&\vdots&\ddots\\
& & &\frac{1-a_0^N}{1-a_0} & &-\frac{1-a_{S -1}^N}{1-a_{S -1}} \\
&& && &&&&-\frac{1-b_S^D}{1-b_S}\\
1-b_1\\
&1-b_2\\
&&\ddots \\
&&&0&&&\textbf{B}\\
&&&&\ddots\\
&&&&&1-b_{S -1}
\end{vmatrix}\\
=&(-1)^{j+S+1}\dfrac{1-b_S^D}{1-b_S}\begin{vmatrix}
1-a_1^Nb_1^D&& &-\frac{(1-a_j)(1-a_0^N)}{1-a_0}  \\
 & 1-a_2^Nb_2^D &&-\frac{(1-a_j)(1-a_0^N)}{1-a_0}   \\
 &&\ddots  &\vdots \\
 &&&-\frac{(1-a_j)(1-a_0^N)}{1-a_0}   \\
 &&&\vdots&\ddots\\
& & &-\frac{(1-a_j)(1-a_0^N)}{1-a_0}  & &1-a_{S-1}^Nb_{S-1}^D \\
\end{vmatrix}\\
=&(-1)^{j+S}\dfrac{1-b_S^D}{1-b_S}\dfrac{(1-a_j)(1-a_0^N)}{1-a_0}\displaystyle\prod_{i=1, 2, ..., \hat{j}, ..., S-1}(1-a_i^Nb_i^D)\\
=&(-1)^{j+S}\dfrac{(1-a_j)(1-a_0^N)^S}{(1-b_S)(1-a_0)}
\end{split}
\end{equation*}
where $\textbf{A} = \text{diag}\left( -a_1^N\frac{1-b_1^D}{1-b_1},  -a_2^N\frac{1-b_2^D}{1-b_2}, ...,  -a_{S-1}^N\frac{1-b_{S-1}^D}{1-b_1}\right)$ and $\textbf{B} = \text{diag}(a_1-1, a_2-1, ..., a_{S-1}-1)$ are diagonal matrix blocks.

$\textbf{Type III.j:}$ Here the $(S+j)^\text{th}$ column is deleted, where $1\leq j \leq S-1$. We exchange columns so that the $(S+j-1)^\text{st}$ column in the Type I determinant is replaced by the first column of the Alexander matrix.

\begin{equation*}
\begin{split}
&(-1)^{S+j-1}\begin{vmatrix}
& -\frac{1-b_1^D}{1-b_1} && &\frac{1-a_0^N}{1-a_0}  \\
 && -\frac{1-b_2^D}{1-b_2} &&\frac{1-a_0^N}{1-a_0} \\
 &&&\ddots  &\vdots \\
\textbf{C}& &&&\frac{1-a_0^N}{1-a_0}\\
& &&&\vdots&\ddots\\
&& & &\frac{1-a_0^N}{1-a_0} & &-\frac{1-b_{S -1}^D}{1-b_{S -1}} \\
&&& & &&&-\frac{1-b_S^D}{1-b_S}\\
&a_1-1\\
&&a_2-1\\
&&&\ddots \\
\textbf{D}&&&&0\\
&&&&&\ddots\\
&&&&&&a_{S -1}-1
\end{vmatrix}\\
=&(-1)^j\dfrac{1-b_S^D}{1-b_S}\begin{vmatrix}
1-a_1^Nb_1^D&& &-\frac{(1-b_j)(1-a_0^N)}{1-a_0}  \\
 & 1-a_2^Nb_2^D &&-\frac{(1-b_j)(1-a_0^N)}{1-a_0}   \\
 &&\ddots  &\vdots \\
 &&&-\frac{(1-b_j)(1-a_0^N)}{1-a_0}   \\
 &&&\vdots&\ddots\\
& & &-\frac{(1-b_j)(1-a_0^N)}{1-a_0}  & &1-a_{S-1}^Nb_{S-1}^D \\
\end{vmatrix}\\
=&(-1)^{j+1}\dfrac{1-b_S^D}{1-b_S}\dfrac{(1-b_j)(1-a_0^N)}{1-a_0}\displaystyle\prod_{i=1, 2, ..., \hat{j}, ..., S-1}(1-a_i^Nb_i^D)\\
=&(-1)^{j}\dfrac{(1-b_j)(1-a_0^N)^S}{(1-b_S)(1-a_0)}
\end{split}
\end{equation*}
where $\textbf{C} = \text{diag}\left( \frac{1-a_1^N}{1-a_1},  \frac{1-a_2^N}{1-a_2}, ..., \frac{1 -a_{S-1}^N}{1-a_1}\right)$ and $\textbf{D} = \text{diag}(1-b_1, 1-b_2-1, ..., 1-b_{S-1})$ are diagonal matrix blocks.

$\textbf{Type IV:}$ Here the $2S^\text{th}$ column is deleted.
\begin{equation*}
\begin{split}
&\begin{vmatrix}
\frac{1-a_0^N}{1-a_0} & -\frac{1-a_1^N}{1-a_1} &&&& -a_1^N\frac{1-b_1^D}{1-b_1} \\
\frac{1-a_0^N}{1-a_0} & & -\frac{1-a_2^N}{1-a_2} &&&&-a_2^N\frac{1-b_2^D}{1-b_2} \\
\vdots & &&\ddots  &&&&\ddots\\
\frac{1-a_0^N}{1-a_0} & && &-\frac{1-a_{S -1}^N}{1-a_{S -1}} &&&&-a_{S -1}^N\frac{1-b_{S -1}^D}{1-b_{S -1}}\\
\frac{1-a_0^N}{1-a_0} \\
& 1-b_1&&&&a_1-1\\
&&1-b_2&&&&a_2-1 \\
&&&\ddots &&&&\ddots \\
&&&&1-b_{S -1}&&&&a_{S -1} -1
\end{vmatrix}\\
=&(-1)^{S+1}\dfrac{1-a_0^N}{1-a_0}\displaystyle\prod_{i=1}^{S-1}(1-a_i^Nb_i^D)\\
=&(-1)^{S+1}\dfrac{(1-a_0^N)^S}{1-a_0}\
\end{split}
\end{equation*}

The Alexander ideal is thus given by
$$\left( \dfrac{(1-a_0^N)^S}{1-b_S}, \left\{\dfrac{(1-a_j)(1-a_0^N)^S}{(1-b_S)(1-a_0)}\right\}_{j=1}^{S-1},\left\{\dfrac{(1-b_j)(1-a_0^N)^S}{(1-b_S)(1-a_0)}\right\}_{j=1}^{S-1}, \dfrac{(1-a_0^N)^S}{1-a_0}\right)$$
The abelianization map $\phi$ is given by the linking numbers, namely,
\begin{equation*}
\begin{split}
\text{lk}(a_0, L) &=SD,\\
\text{lk}(a_j, L) &=(S-j)D,\\
\text{lk}(b_j, L) &= jN
\end{split}
\end{equation*} and $\phi(x) = t^{\text{lk}(x, L)}$.

The Alexander ideal under the abelianization map is given by
\begin{equation*}
\begin{split}
 &\left( \dfrac{(1-t^{SND})^S}{1-t^{SN}}, \left\{\dfrac{(1-t^{(S-j)D})(1-t^{SND})^S}{(1-t^{SN})(1-t^{SD})}\right\}_{j=1}^{S-1},\left\{\dfrac{(1-t^{jN})(1-t^{SND})^S}{(1-t^{SN})(1-t^{SD})}\right\}_{j=1}^{S-1}, \dfrac{(1-t^{SND})^S}{1-t^{SD}}\right)\\
=&\dfrac{(1-t^{SND})^S}{(1-t^{SN})(1-t^{SD})}\left(1-t^{SD}, \left\{1-t^{(S-j)D}\right\}_{j=1}^{S-1},\left\{1-t^{jN}\right\}_{j=1}^{S-1}, 1-t^{SN}\right)\\
=&\dfrac{(1-t^{SND})^S}{(1-t^{SN})(1-t^{SD})}\left(1-t^{SD}, 1-t^D, \left\{1-t^{(S-j)D}\right\}_{j=2}^{S-1}, 1-t^N, \left\{1-t^{jN}\right\}_{j=2}^{S-1}, 1-t^{SN}\right)\\
=&\dfrac{(1-t^{SND})^S}{(1-t^{SN})(1-t^{SD})}(1-t)
\end{split}
\end{equation*}since $N$ and $D$ are coprime integers.

Then, to obtain the Alexander polynomial, we have to multiply by a power of $t$ so that the function gives the same Laurent polynomial in $\sqrt t$ when evaluated in $t$ and $t^{-1}$: the power is $-\frac12(1+S^2ND - SN - SD) = -\frac 12 (SN-1)(SD-1)$. Besides it is customary to add the factor $(-1)^{S-1}$ so that the leading term has positive coefficient.
 
Hence we have
$$\Delta_L (t) = (-1)^{S-1}t^{-\frac{(SN-1)(SD-1)}2} \dfrac{(1-t)(1-t^{SND})^S}{(1-t^{SN})(1-t^{SD})}$$
Now we note that $K_{-SN,SD}$ is the mirror image of $K_{SN,SD}$ \cite{Murasugi}, and mirror images have the same Alexander polynomial \cite{Collins}, so we have the desired result.

\subsection{Tst links $[\Phi^\mu(n(\tau), d(\tau), M)]$ with odd $M$}
Next we consider the Alexander polynomials of the tst link $[\Phi^\mu(n(\tau),d(\tau), M)]$ when $M$ is odd, which we conjectured in \cite{HLMA} to be
$$\Delta_{[\Phi^\mu(n(\tau),d(\tau), M)]} (t)= (-1)^{\mu P}t^{-\frac{(\mu PN(\tau)-1)(\mu P D(\tau)+1)+1}2} \dfrac{(1-t)(1-t^{(\mu PD(\tau) +1)N(\tau)})^{\mu P}}{1-t^{\mu PD(\tau) +1}}$$
Here we add a factor of $(-1)^{\mu P}$. We shall now prove this result.

When $M$ is odd, the link is given by a $K_{-\mu PN(\tau), \mu PD(\tau)}$ on the surface of the torus and inside the solid torus bounded by it there is a trivial knot along the longitudinal direction. As above we shall write $N=N(\tau)$, $D=D(\tau)$, $\mu P = S$ and consider the mirror image $L$ where $K_{SN, SD}$ is placed on the surface of the torus instead of $K_{-SN, SD}$.

Like the case where $M$ is even, we introduce $S$ solid tori $U_1, U_2, ..., U_S$ with the same longitudinal circle where $U_S \subset \cdots \subset U_2 \subset U_1$ and place a component $K_{N,D}$ of the torus link on each of the boundaries $\partial U_i$. By Seifert-van Kampen Theorem the fundamental group of $S^3\backslash L$ is given by
$$\pi_1(S^3\backslash L) \cong \left\langle
\begin{split}
a_0, a_1, a_2, ..., a_S,\\
b_1, b_2, ..., b_S
\end{split}
\biggr|
\begin{split}
a_0^N = a_i^Nb_i^D, & 1\leq i \leq S\\
a_jb_j=b_ja_j, & 1\leq j \leq S
\end{split}
\right\rangle $$
Now we set the relators $r_i = a_0^N (a_i^Nb_i^D)^{-1}$ and $s_j = (a_jb_j)(b_ja_j)^{-1}$. The Fox derivatives are given by
\begin{equation*}
\begin{split}
\dfrac{\partial r_i}{\partial a_0} &= \dfrac{1-a_0^N}{1-a_0}\\
\dfrac{\partial r_i}{\partial a_i} &= -\dfrac{1-a_i^N}{1-a_i}\\
\dfrac{\partial r_i}{\partial b_i} &= -a_i^N\dfrac{1-b_i^D}{1-b_i}\\
\dfrac{\partial s_j}{\partial a_j} &= 1-b_j \\
\dfrac{\partial s_j}{\partial b_j} &= a_j - 1
\end{split}
\end{equation*}
where $1\leq i,j \leq S$ and other Fox derivatives are zero.

The Alexander matrix is given by
$$\begin{bmatrix}
\frac{1-a_0^N}{1-a_0} & -\frac{1-a_1^N}{1-a_1} &&&& -a_1^N\frac{1-b_1^D}{1-b_1} \\
\frac{1-a_0^N}{1-a_0} & & -\frac{1-a_2^N}{1-a_2} &&&&-a_2^N\frac{1-b_2^D}{1-b_2} \\
\vdots & &&\ddots  &&&&\ddots\\
\frac{1-a_0^N}{1-a_0} & && &-\frac{1-a_S^N}{1-a_S} &&&&-a_S^N\frac{1-b_S^D}{1-b_S}\\
& 1-b_1&&&&a_1-1\\
&&1-b_2&&&&a_2-1 \\
&&&\ddots &&&&\ddots \\
&&&&1-b_S&&&&a_S -1
\end{bmatrix}$$
The Alexander ideal is generated by the following types of determinants:

$\textbf{Type I:}$ Here the first column is deleted.
\begin{equation*}
\begin{split}
&\begin{vmatrix}
 -\frac{1-a_1^N}{1-a_1} &&&& -a_1^N\frac{1-b_1^D}{1-b_1} \\
 & -\frac{1-a_2^N}{1-a_2} &&&&-a_2^N\frac{1-b_2^D}{1-b_2} \\
 &&\ddots  &&&&\ddots\\
 && &-\frac{1-a_S^N}{1-a_S} &&&&-a_S^N\frac{1-b_S^D}{1-b_S}\\
1-b_1&&&&a_1-1\\
&1-b_2&&&&a_2-1 \\
&&\ddots &&&&\ddots \\
&&&1-b_S&&&&a_S -1
\end{vmatrix}\\
=& \prod_{i=1}^S(1-a_i^Nb_i^D)\\
=& (1-a_0^N)^S
\end{split}
\end{equation*}
\clearpage
$\textbf{Type II.j:}$ Here the $(j+1)^\text{st}$ column is deleted, where $1\leq j \leq S$. We exchange columns so that the $j^\text{th}$ column in the Type I determinant is replaced by the first column of the Alexander matrix.

\begin{equation*}
\begin{split}
&(-1)^j\begin{vmatrix}
 -\frac{1-a_1^N}{1-a_1} && &\frac{1-a_0^N}{1-a_0}  \\
 & -\frac{1-a_2^N}{1-a_2} &&\frac{1-a_0^N}{1-a_0} \\
 &&\ddots  &\vdots \\
 &&&\frac{1-a_0^N}{1-a_0}&&&\textbf{A} \\
 &&&\vdots&\ddots\\
& & &\frac{1-a_0^N}{1-a_0} & &-\frac{1-a_S^N}{1-a_S} \\
1-b_1\\
&1-b_2\\
&&\ddots \\
&&&0&&&\textbf{B}\\
&&&&\ddots\\
&&&&&1-b_S
\end{vmatrix}\\
=&(-1)^j\begin{vmatrix}
1-a_1^Nb_1^D&& &-\frac{(1-a_j)(1-a_0^N)}{1-a_0}  \\
 & 1-a_2^Nb_2^D &&-\frac{(1-a_j)(1-a_0^N)}{1-a_0}   \\
 &&\ddots  &\vdots \\
 &&&-\frac{(1-a_j)(1-a_0^N)}{1-a_0}   \\
 &&&\vdots&\ddots\\
& & &-\frac{(1-a_j)(1-a_0^N)}{1-a_0}  & &1-a_S^Nb_S^D \\
\end{vmatrix}\\
=&(-1)^{j+1}\dfrac{(1-a_j)(1-a_0^N)}{1-a_0}\displaystyle\prod_{i=1, 2, ..., \hat{j}, ..., S}(1-a_i^Nb_i^D)\\
=&(-1)^{j+1}\dfrac{(1-a_j)(1-a_0^N)^S}{1-a_0}
\end{split}
\end{equation*}
where $\textbf{A} = \text{diag}\left( -a_1^N\frac{1-b_1^D}{1-b_1},  -a_2^N\frac{1-b_2^D}{1-b_2}, ...,  -a_S^N\frac{1-b_S^D}{1-b_1}\right)$ and $\textbf{B} = \text{diag}(a_1-1, a_2-1, ..., a_S-1)$ are diagonal matrix blocks.
\clearpage
$\textbf{Type III.j:}$ Here the $(S+j-1)^\text{st}$ column is deleted, where $1\leq j \leq S$. We exchange columns so that the $(S+j)^\text{th}$ column in the Type I determinant is replaced by the first column of the Alexander matrix.

\begin{equation*}
\begin{split}
&(-1)^{S+j}\begin{vmatrix}
& -\frac{1-b_1^D}{1-b_1} && &\frac{1-a_0^N}{1-a_0}  \\
 && -\frac{1-b_2^D}{1-b_2} &&\frac{1-a_0^N}{1-a_0} \\
 &&&\ddots  &\vdots \\
\textbf{C}& &&&\frac{1-a_0^N}{1-a_0}\\
& &&&\vdots&\ddots\\
&& & &\frac{1-a_0^N}{1-a_0} & &-\frac{1-b_S^D}{1-b_S} \\
&a_1-1\\
&&a_2-1\\
&&&\ddots \\
\textbf{D}&&&&0\\
&&&&&\ddots\\
&&&&&&a_S-1
\end{vmatrix}\\
=&(-1)^{S+j}\begin{vmatrix}
1-a_1^Nb_1^D&& &-\frac{(1-b_j)(1-a_0^N)}{1-a_0}  \\
 & 1-a_2^Nb_2^D &&-\frac{(1-b_j)(1-a_0^N)}{1-a_0}   \\
 &&\ddots  &\vdots \\
 &&&-\frac{(1-b_j)(1-a_0^N)}{1-a_0}   \\
 &&&\vdots&\ddots\\
& & &-\frac{(1-b_j)(1-a_0^N)}{1-a_0}  & &1-a_S^Nb_S^D \\
\end{vmatrix}\\
=&(-1)^{S+j+1}\dfrac{(1-b_j)(1-a_0^N)}{1-a_0}\displaystyle\prod_{i=1, 2, ..., \hat{j}, ..., S}(1-a_i^Nb_i^D)\\
=&(-1)^{S+j+1}\dfrac{(1-b_j)(1-a_0^N)^S}{1-a_0}
\end{split}
\end{equation*}
where $\textbf{C} = \text{diag}\left( \frac{1-a_1^N}{1-a_1},  \frac{1-a_2^N}{1-a_2}, ..., \frac{1 -a_S^N}{1-a_1}\right)$ and $\textbf{D} = \text{diag}(1-b_1, 1-b_2-1, ..., 1-b_S)$ are diagonal matrix blocks.

Therefore the Alexander ideal is
$$\left( (1-a_0^N)^S, \left\{\dfrac{(1-a_j)(1-a_0^N)^S}{1-a_0}\right\}_{j=1}^S, \left\{\dfrac{(1-b_j)(1-a_0^N)^S}{1-a_0}\right\}_{j=1}^S \right)$$

The abelianization map is given by
$\phi(x)=t^{\text{lk}(x, L)}$:
\begin{equation*}
\begin{split}
\phi(a_0)&=t^{SD+1}\\
\phi(a_j)&= t^{(S-j)D+1}\\
\phi(b_j)&= t^{jN}
\end{split}
\end{equation*}
Under this map the Alexander ideal is
\begin{equation*}
\begin{split}
&\left( (1-t^{(SD+1)N})^S, \left\{\dfrac{(1-t^{(S-j)D+1})(1-t^{(SD+1)N})^S}{1-t^{SD+1}}\right\}_{j=1}^S, \left\{\dfrac{(1-t^{jN})(1-t^{(SD+1)N})^S}{1-t^{SD+1}}\right\}_{j=1}^S \right)\\
=&\dfrac{(1-t^{(SD+1)N})^S}{1-t^{SD+1}}\left( 1-t^{SD+1}, \left\{1-t^{(S-j)D+1}\right\}_{j=1}^S, \left\{1-t^{jN}\right\}_{j=1}^S \right)\\
=&\dfrac{(1-t^{(SD+1)N})^S}{1-t^{SD+1}}\left( 1-t^{SD+1}, \left\{1-t^{(S-j)D+1}\right\}_{j=1}^{S-1}, 1-t, \left\{1-t^{jN}\right\}_{j=1}^S \right)\\
=&\dfrac{(1-t^{(SD+1)N})^S}{1-t^{SD+1}}(1-t)
\end{split}
\end{equation*}
The generator of the ideal is thus given by
$$\dfrac{(1-t)(1-t^{(SD+1)N})^S}{1-t^{SD+1}}$$
The power of $t$ to multiply is $t^{-\frac 12(1+(SD+1)NS-(SD+1))} = t^{- \frac {(SD+1)(SN-1)+1}2}$. We also multiply a  power of $-1$ so that the leading term has positive coefficient, which is given by $(-1)^S$. Hence
$$\Delta_L (t)= (-1)^St^{-\frac{(SN-1)(SD+1)+1}2} \dfrac{(1-t)(1-t^{(\mu PD(\tau) +1)N(\tau)})^{\mu P}}{1-t^{\mu PD(\tau) +1}}$$

For odd $M$, the link $[\Phi^\mu (n(\tau), d(\tau), M)]$ is the mirror image of $L$, and since mirror images have the same Alexander polynomial \cite{Collins}, the claimed result follows.

\subsection{Stst links $[\Phi_{j=1}^k(n(\tau_j),d(\tau_j), M_j, \mu_j)]$ with odd $M_k$}
We now turn to more general surgerized twisted solid tori links (stst links) $[\Phi_{j=1}^k (n(\tau_j), d(\tau_j), M_j, \mu_j)]$. We introduce the notations $N(\tau_j)=N_j, D(\tau_j)=D_j, \mu_j P_j=S_j$.  We first consider the case where $M_k$ is odd. We use the fact that mirror images have the same Alexander polynomial \cite{Collins}, and consider $\sum_{j=1}^k S_j$ solid tori with the same longitudinal circle, $U_{k, S_k} \subset \cdots \subset U_{k, 1} \subset  \cdots \subset U_{1, S_1}\subset \cdots \subset U_{1, 1}$, inside $S^3$. Then we place the $S_1$ components of the torus link $K_{S_1N_1,S_1D_1}$ on the surfaces $\partial U_{S_1, 1}, ..., \partial U_{S_1, S_1}$, and so on, up to the $S_k$ components of the torus link $K_{S_kN_k,S_kD_k}$ on the surfaces $\partial U_{S_k, 1}, ..., \partial U_{S_k, S_k}$, and place a trivial knot inside the solid torus $U_{S_k,S_k}$. We denote the collection of all these torus knots and the trivial knot by $L$. The fundamental group of the complement of $L$ is given by
$$\pi_1(S^3\backslash L)\cong \left \langle 
\begin{split}
a_0, a_{1,1}, \cdots , a_{1,S_1}, \cdots ,a_{k,1}, \cdots, a_{k, S_k},\\
b_{1,1}, \cdots, b_{1,S_1}, \cdots, b_{k,1}, \cdots, b_{k, S_k}
\end{split} \Biggr |
\begin{split}
 a_0^{N_1} = a_{1,i_1}^{N_1}b_{1,i_1}^{D_1},   a_{l-1, S_{l-1 }}^{N_l}b_{l-1, S_{l-1}}^{D_l} = a_{l,i_l}^{N_l}b_{l,i_l}^{D_l}\\
 a_{m,n}b_{m,n} = b_{m,n}a_{m,n}
\end{split}\right \rangle$$
where $1\leq i_j \leq S_j$, $2\leq l \leq k$, $1\leq m \leq k$, and $1\leq n \leq S_m$.

The Alexander matrix is given by
$$
\begin{bmatrix}
\textbf{v}_0&A_1 &          &               &        &C_1\\
                     &B_1 &A_2    &               &        &D_1 & C_2 \\
                     &       &\ddots&\ddots &        &        &\ddots&\ddots\\
                     &       &             & B_{k-1}&A_k&        && D_{k-1}&C_k \\
       \textbf{0}              &      &E          &               &       &        & F
\end{bmatrix}
$$
where $\textbf{v}_0 = \frac{1-a_0^N}{1-a_0}\begin{bmatrix}
1&1&\cdots &1
\end{bmatrix}^T$ has $S_1$ rows,
\begin{equation*}
\begin{split}
A_j &= -\text{diag}\left( \dfrac{1-a_{j,1}^{N_j}}{1-a_{j,1}},\dfrac{1-a_{j,2}^{N_j}}{1-a_{j,2}}, \cdots, \dfrac{1-a_{j,S_j}^{N_j}}{1-a_{j,S_j}}\right)\\
B_j &= \dfrac{1-a_{j, S_j}^{N_{j+1}}}{1-a_{j, S_j}}\begin{bmatrix}
0 &\cdots & 0 &  1\\
\vdots&\ddots&\vdots  &\vdots\\
0 &\cdots & 0 &1
\end{bmatrix}\\
C_j &= -\text{diag}\left(a_{j,1}^{N_j} \dfrac{1-b_{j,1}^{D_j}}{1-b_{j,1}},a_{j,2}^{N_j} \dfrac{1-b_{j,2}^{D_j}}{1-b_{j,2}}, \cdots, a_{j, S_j}^{N_j} \dfrac{1-b_{j,S_j}^{D_j}}{1-b_{j,S_j}}\right)\\
D_j &= a_{j, S_j}^{N_{j+1}} \dfrac{1-b_{j, S_j}^{D_{j+1}}}
{1-b_{j, S_j}}
\begin{bmatrix}
0 &\cdots & 0 &  1\\
\vdots&\ddots&\vdots  &\vdots\\
0 &\cdots & 0 &1
\end{bmatrix}\\
E &= \text{diag}(1-b_{1,1}, \cdots, 1-b_{1, S_1},\cdots, 1-b_{k, 1},\cdots, 1-b_{k, S_k})\\
F &= \text{diag}(a_{1,1}-1, \cdots, a_{1, S_1}-1,\cdots, a_{k, 1}-1,\cdots, a_{k, S_k}-1)\\
\end{split}
\end{equation*}
and $B_j$ and $D_j$ both have $S_{j+1}$ rows and $S_j$ columns.
The Alexander ideal is generated by the following types of determinants:

$\textbf{Type I:}$ The first column of the Alexander matrix is deleted. The determinant is given by
$$\prod_{j=1}^k \prod_{i=1}^{S_j} \left(1-a_{j, i}^{N_j} b_{j, i}^{D_j}\right) =\prod_{j=1}^k \left(1-a_{j, 1}^{N_j} b_{j, 1}^{D_j}\right)^{S_j} $$

$\textbf{Type II.i:}$ The $(i+1)^\text{st}$ column of the Alexander matrix is deleted, where $1\leq i \leq S_1$. The determinant is given by
\begin{equation*}
\begin{split}
&(-1)^{i+1}\left(\prod_{j=2}^k \prod_{m=1}^{S_j} \left(1-a_{j, m}^{N_j} b_{j, m}^{D_j}\right)\right) \left(\prod_{l = 1, 2, ..., \hat{i}, ..., S_1}
\left(1-a_{1, i}^{N_1} b_{1, i}^{D_1}\right) \right)\dfrac{(1-a_{1,i})(1-a_0^{N_1})}{1-a_0}
\\
=&(-1)^{i+1}\dfrac{(1-a_{1,i})\displaystyle\prod_{j=1}^k \left(1-a_{j, 1}^{N_j} b_{j, 1}^{D_j}\right)^{S_j} }{1-a_0}
\end{split}
\end{equation*}
The determinant of the matrix where the $(i+1)^\text{st}$ column of the Alexander matrix is deleted is zero for $S_1+1 \leq i \leq 1+\sum_{j=1}^k S_k$.

$\textbf{Type III.i:}$ The $(i+1+\sum_{j=1}^k S_k)^\text{th}$ column of the Alexander matrix is deleted, where $1\leq i \leq S_1$. The determinant is given by
\begin{equation*}
\begin{split}
&(-1)^{i+1+\sum_{j=1}^k S_k}\left(\prod_{j=2}^k \prod_{m=1}^{S_j} \left(1-a_{j, m}^{N_j} b_{j, m}^{D_j}\right)\right) \left(\prod_{l = 1, 2, ..., \hat{i}, ..., S_1}
\left(1-a_{1, i}^{N_1} b_{1, i}^{D_1}\right) \right)\dfrac{(1-b_{1,i})(1-a_0^{N_1})}{1-a_0}
\\
=&(-1)^{i+1+\sum_{j=1}^k S_k}\dfrac{(1-b_{1,i})\displaystyle\prod_{j=1}^k \left(1-a_{j, 1}^{N_j} b_{j, 1}^{D_j}\right)^{S_j} }{1-a_0}
\end{split}
\end{equation*}
The determinant of the matrix where the $(i+1+\sum_{j=1}^k S_k)^\text{th}$ column of the Alexander matrix is deleted  is zero for $S_1+1 \leq i \leq \sum_{j=1}^k S_k$.

The abelianization map $\phi$ is given by $\phi(a_0) = t^{1+\sum_{j=1}^k S_j D_j}$, $\phi(a_{j, i}) = t^{1+(S_j-i)D_j + \sum_{l=j+1}^k S_lD_l}$ and $\phi(b_{j,i}) = t^{iN_j +\sum_{l=1}^{j-1}S_lN_l}$, where the sum $\sum_{l=k+1}^k$ and $\sum_{l=1}^0$ are interpreted as $0$.

Writing
\begin{equation*}
\begin{split}
\Pi &= \prod_{j=1}^k \left(1-t^{\left(1+(S_j-1)D_j + \sum_{l=j+1}^k S_lD_l\right)N_j +  \left(N_j +\sum_{l=1}^{j-1}S_lN_l\right)D_j} \right)^{S_j} \\
&= \prod_{j=1}^k \left(1-t^{N_j+S_jN_jD_j + N_j\sum_{l=j+1}^k S_lD_l +  D_j\sum_{l=1}^{j-1}S_lN_l} \right)^{S_j}
\end{split}
\end{equation*}
the Alexander ideal under the abelianization map is given by

\begin{equation*}
\begin{split}
&\left(\Pi, 
\left\{\dfrac{\left(1-t^{1-iD_1 + \sum_{l=1}^k S_lD_l}\right)\Pi }{1-t^{1+\sum_{j=1}^k S_j D_j}}\right\}_{i=1}^{S_1}, 
\left\{\dfrac{(1-t^{iN_1})\Pi}{1-t^{1+\sum_{j=1}^k S_j D_j}}\right\}_{i=1}^{S_1}
\right)\\
=&\dfrac\Pi{1-t^{1+\sum_{j=1}^k S_j D_j}}\left(1-t^{1+\sum_{j=1}^k S_j D_j}, 
\left\{1-t^{1+(S_1-i)D_1 + \sum_{l=2}^k S_lD_l}\right\}_{i=1}^{S_1}, 
\left\{1-t^{iN_1}\right\}_{i=1}^{S_1}
\right)\\
=&\dfrac\Pi{\left[1+\sum_{j=1}^k S_j D_j\right]}\left(\left[1+\sum_{j=1}^k S_j D_j\right], 
\left\{\left[1+(S_1-i)D_1 + \sum_{l=1}^k S_lD_l\right]\right\}_{i=1}^{S_1}, 
\left\{[iN_1]\right\}_{i=1}^{S_1}
\right)
\end{split}
\end{equation*}
where to write more compactly, we introduce the notation $[\lambda] = 1-t^\lambda$ for integers $\lambda$. 

We note that the following diagram is commutative:

$$\begin{tikzcd}
\mathbb{Z}^n \arrow{r}{[\cdot]^n} \arrow[swap]{d}{\text{gcd}} & {\mathbb{Z}[t]^n} \arrow{d}{\text{gcd}} \\%
\mathbb{Z}\arrow{r}{[\cdot]}& {\mathbb{Z}[t]}
\end{tikzcd}$$We thus have
\begin{equation*}
\begin{split}
&\text{gcd}\left(\left[1+\sum_{j=1}^k S_j D_j\right], 
\left\{\left[1+(S_1-i)D_1 + \sum_{l=1}^k S_lD_l\right]\right\}_{i=1}^{S_1}, 
\left\{[iN_1]\right\}_{i=1}^{S_1}
\right) \\
=& \left[\text{gcd}\left(1+\sum_{j=1}^k S_j D_j, 
\left\{1+(S_1-i)D_1 + \sum_{l=1}^k S_lD_l\right\}_{i=1}^{S_1}, 
\left\{iN_1\right\}_{i=1}^{S_1}\right)\right]\\
=& \left[\text{gcd}\left(
\text{gcd}\left(1+(S_1-1)D_1 + \sum_{l=1}^k S_lD_l, 1+(S_1-2)D_1 + \sum_{l=1}^k S_lD_l\right),
N_1, \cdots\right)\right]\\
=& \left[\text{gcd}\left(
\text{gcd}\left(D_1 , 1+(S_1-2)D_1 + \sum_{l=1}^k S_lD_l\right),
N_1, \cdots\right)\right]\\
=& \left[\text{gcd}\left(D_1,N_1, \cdots\right)\right]\\
= &[1] \\
=&1-t
\end{split}
\end{equation*}where we assumed $S_1\geq 2$.

The generator of the Alexander ideal is thus
$$\dfrac{(1-t)\displaystyle\prod_{j=1}^k \left(1-t^{N_j+S_jN_jD_j + N_j\sum_{l=j+1}^k S_lD_l +  D_j\sum_{l=1}^{j-1}S_lN_l} \right)^{S_j}} {1-t^{1+\sum_{j=1}^k S_j D_j}} $$
Multiplying by powers of $t$ and $(-1)$, the Alexander polynomial of $L$ is
$$\Delta_L(t) = (-1)^{\sum_{j=1}^k S_j } 
t^\frac \sigma 2
\dfrac{(1-t)\displaystyle\prod_{j=1}^k \left(1-t^{N_j+S_jN_jD_j + N_j\sum_{l=j+1}^k S_lD_l +  D_j\sum_{l=1}^{j-1}S_lN_l} \right)^{S_j}} {1-t^{1+\sum_{j=1}^k S_j D_j}}$$
where $\sigma = \sum_{j=1}^k\left( S_j D_j - (N_j+S_jN_jD_j + N_j\sum_{l=j+1}^k S_lD_l +  D_j\sum_{l=1}^{j-1}S_lN_l)\right)$.

\subsection{Stst links $[\Phi_{j=1}^k(n(\tau_j),d(\tau_j), M_j, \mu_j)]$ with even $M_k$}
We now consider stst links $[\Phi_{j=1}^k (n(\tau_j), d(\tau_j), M_j, \mu_j)]$ with even $M_k$. Writing $N(\tau_j)=N_j, D(\tau_j)=D_j, \mu_j P_j=S_j$, we also use the fact that mirror images have the same Alexander polynomial \cite{Collins}.  The link $[\Phi_{j=1}^k (n(\tau_j), d(\tau_j), M_j, \mu_j)]$  contains $\sum_{j=1}^k S_j$ torus knots. We denote the collection of the mirror images of these torus knots by $L$. Using the technique as above, the fundamental group of the complement of $L$ is given by
$$\pi_1(S^3\backslash L)\cong \left \langle 
\begin{split}
a_0, a_{1,1}, \cdots , a_{k-1,S_{k-1}}, a_{k, 1},\cdots, a_{k, S_k-1},\\
b_{1,1}, \cdots, b_{1,S_1}, \cdots, b_{k,1}, \cdots, b_{k, S_k}
\end{split} \Biggr |
\begin{split}
 a_0^{N_1} = a_{1,i_1}^{N_1}b_{1,i_1}^{D_1},   a_{l-1, S_{l-1 }}^{N_l}b_{l-1, S_{l-1}}^{D_l} = a_{l,i_l}^{N_l}b_{l,i_l}^{D_l}\\
a_{k-1, S_{k-1}}^{N_{k-1}}b_{k-1, S_{k-1}}^{D_{k-1}} = b_{k,S_k}^{D_k}, a_{m,n}b_{m,n} = b_{m,n}a_{m,n}
\end{split}\right \rangle$$
where $1\leq i_j \leq S_j$ for $1\leq j \leq k-1$ and $1\leq i_k \leq S_k-1$, $2\leq l \leq k$, $1\leq n \leq S_m$ for $1\leq m \leq k-1$, and $1\leq n \leq S_k-1$ for $m=k$.

The Alexander matrix is given by
$$
\begin{bmatrix}
\textbf{v}_0&A_1 &          &               &        &C_1\\
                     &B_1 &A_2    &               &        &D_1 & C_2 \\
                     &       &\ddots&\ddots &        &        &\ddots&\ddots\\
                     &       &             & B_{k-1}&A_k&        && D_{k-1}&C_k \\
       \textbf{0}              &       &E          &        &       &        & F
\end{bmatrix}
$$
where $\textbf{v}_0 = \frac{1-a_0^N}{1-a_0}\begin{bmatrix}
1&1&\cdots &1
\end{bmatrix}^T$ has $S_1$ rows,
\begin{equation*}
\begin{split}
A_j &= -\text{diag}\left( \dfrac{1-a_{j,1}^{N_j}}{1-a_{j,1}},\dfrac{1-a_{j,2}^{N_j}}{1-a_{j,2}}, \cdots, \dfrac{1-a_{j,S_j}^{N_j}}{1-a_{j,S_j}}\right)\\
A_k&= \begin{bmatrix}
-\text{diag}\left( \dfrac{1-a_{k,1}^{N_k}}{1-a_{k,1}},\dfrac{1-a_{k,2}^{N_k}}{1-a_{k,2}}, \cdots, \dfrac{1-a_{k,S_k-1}^{N_k}}{1-a_{j,S_k-1}}\right)\\
\textbf{0}
\end{bmatrix}\\
B_l &= \dfrac{1-a_{j, S_j}^{N_{j+1}}}{1-a_{j, S_j}}\begin{bmatrix}
0 &\cdots & 0 &  1\\
\vdots&\ddots&\vdots  &\vdots\\
0 &\cdots & 0 &1
\end{bmatrix}\\
C_j &= -\text{diag}\left(a_{j,1}^{N_j} \dfrac{1-b_{j,1}^{D_j}}{1-b_{j,1}},a_{j,2}^{N_j} \dfrac{1-b_{j,2}^{D_j}}{1-b_{j,2}}, \cdots, a_{j, S_j}^{N_j} \dfrac{1-b_{j,S_j}^{D_j}}{1-b_{j,S_j}}\right)\\
C_k &= -\text{diag}\left(a_{k,1}^{N_k} \dfrac{1-b_{k,1}^{D_k}}{1-b_{k,1}},a_{k,2}^{N_k} \dfrac{1-b_{k,2}^{D_k}}{1-b_{k,2}}, \cdots, a_{k, S_k-1}^{N_k-1} \dfrac{1-b_{k,S_k-1}^{D_k}}{1-b_{k,S_k-1}}, \dfrac{1-b_{k,S_k}^{D_k}}{1-b_{k,S_k}}\right)\\
D_l &= a_{j, S_j}^{N_{j+1}} \dfrac{1-b_{j, S_j}^{D_{j+1}}}
{1-b_{j, S_j}}
\begin{bmatrix}
0 &\cdots & 0 &  1\\
\vdots&\ddots&\vdots  &\vdots\\
0 &\cdots & 0 &1
\end{bmatrix}\\
E &= \text{diag}(1-b_{1,1}, \cdots, 1-b_{1, S_1},\cdots, 1-b_{k, 1},\cdots, 1-b_{k, S_k})\\
F &= \text{diag}(a_{1,1}-1, \cdots, a_{1, S_1}-1,\cdots, a_{k, 1}-1,\cdots, a_{k, S_k}-1)\\
\end{split}
\end{equation*}
where $B_l$ and $D_l$ both have $S_{l+1}$ rows and $S_l$ columns, $1\leq l \leq k$ and $1\leq j \leq k-1$ .

The Alexander ideal is generated by the following types of determinants:

$\textbf{Type I:}$ The first column of the Alexander matrix is deleted. The determinant, up to sign, is given by
$$\left(\prod_{j=1}^{k-1} \prod_{i=1}^{S_j} \left(1-a_{j, i}^{N_j} b_{j, i}^{D_j}\right) \right)\dfrac{1-b_{k,S_k}^{D_k}}{1-b_{k,S_k}}\prod_{i=1}^{S_k-1}\left(1-a_{k, i}^{N_k} b_{k, i}^{D_k}\right) 
=\dfrac{\displaystyle\prod_{j=1}^k \left(1-a_{j, 1}^{N_j} b_{j, 1}^{D_j}\right)^{S_j}}{1-b_{k,S_k}} $$

$\textbf{Type II.i:}$ The $(i+1)^\text{st}$ column of the Alexander matrix is deleted, where $1\leq i \leq S_1$. The determinant, up to sign, is given by
\begin{equation*}
\begin{split}
&\left(\prod_{j=2}^k \prod_{m=1}^{S_j} \left(1-a_{j, m}^{N_j} b_{j, m}^{D_j}\right)\right) \left(\prod_{l = 1, 2, ..., \hat{i}, ..., S_1}
\left(1-a_{1, i}^{N_1} b_{1, i}^{D_1}\right) \right)\dfrac{(1-a_{1,i})(1-a_0^{N_1})}{(1-a_0)(1-b_{k,S_k})}
\\
=&\dfrac{(1-a_{1,i})\displaystyle\prod_{j=1}^k \left(1-a_{j, 1}^{N_j} b_{j, 1}^{D_j}\right)^{S_j} }{(1-a_0)(1-b_{k,S_k})}
\end{split}
\end{equation*}
The determinant of the matrix where the $(i+1)^\text{st}$ column of the Alexander matrix is deleted is zero for $S_1+1 \leq i \leq -1+\sum_{j=1}^k S_k$.

$\textbf{Type III.i:}$ The $(i+1+\sum_{j=1}^k S_k)^\text{th}$ column of the Alexander matrix is deleted, where $1\leq i \leq S_1$. Up to sign, the determinant is given by
\begin{equation*}
\begin{split}
&\left(\prod_{j=2}^k \prod_{m=1}^{S_j} \left(1-a_{j, m}^{N_j} b_{j, m}^{D_j}\right)\right) \left(\prod_{l = 1, 2, ..., \hat{i}, ..., S_1}
\left(1-a_{1, i}^{N_1} b_{1, i}^{D_1}\right) \right)\dfrac{(1-b_{1,i})(1-a_0^{N_1})}{(1-a_0)(1-b_{k,S_k})}
\\
=&\dfrac{(1-b_{1,i})\displaystyle\prod_{j=1}^k \left(1-a_{j, 1}^{N_j} b_{j, 1}^{D_j}\right)^{S_j} }{(1-a_0)(1-b_{k,S_k})}
\end{split}
\end{equation*}
The determinant of the matrix where the $(i+\sum_{j=1}^k S_k)^\text{th}$ column of the Alexander matrix is deleted  is zero for $S_1+1 \leq i \leq \sum_{j=1}^k S_k$.

$\textbf{Type IV:}$ The last column of the Alexander matrix is deleted. The determinant of the matrix is zero.

The abelianization map $\phi$ is given by $\phi(a_0) = t^{\sum_{j=1}^k S_j D_j}$, $\phi(a_{j, i}) = t^{(S_j-i)D_j + \sum_{l=j+1}^k S_lD_l}$ and $\phi(b_{j,i}) = t^{iN_j +\sum_{l=1}^{j-1}S_lN_l}$, where the sum $\sum_{l=k+1}^k$ and $\sum_{l=1}^0$ are interpreted as $0$.

Writing 
\begin{equation*}
\begin{split}
\Pi &= \prod_{j=1}^k \left(1- t^{(S_j-1)N_jD_j +N_j \sum_{l=j+1}^k S_lD_l + 
 D_jN_j +D_j\sum_{k=1}^{j-1}S_lN_l}\right)^{S_j} \\
 &= \prod_{j=1}^k \left(1- t^{S_jN_jD_j +N_j \sum_{l=j+1}^k S_lD_l + 
D_j\sum_{k=1}^{j-1}S_lN_l}\right)^{S_j}
\end{split}
\end{equation*}
the Alexander ideal under the abelianization map is given by
\begin{equation*}
\begin{split}
&\left( \dfrac{\Pi}{1-t^{\sum_{j=1}^k S_j N_j}}, 
\left\{\dfrac{\left(1-t^{(S_1-i)D_1 + \sum_{l=2}^k S_lD_l}\right)\Pi }{\left(1-t^{\sum_{j=1}^k S_j D_j}\right)\left(1-t^{\sum_{j=1}^k S_j N_j}\right)}\right\}_{i=1}^{S_1}, 
\left\{\dfrac{\left(1-t^{iN_1}\right)\Pi }{\left(1-t^{\sum_{j=1}^k S_j D_j}\right)\left(1-t^{\sum_{j=1}^k S_j N_j}\right)}\right\}_{i=1}^{S_1}\right)\\
=&\dfrac{\Pi}{\left(1-t^{\sum_{j=1}^k S_j D_j}\right)\left(1-t^{\sum_{j=1}^k S_j N_j}\right)}\left(1-t^{\sum_{j=1}^k S_j D_j},  \{1-t^{(S_1-i)D_1 + \sum_{l=2}^k S_lD_l}\}_{i=1}^{S_1}, \{1-t^{iN_1}\}_{i=1}^{S_1}\right)\\
=&\dfrac{\Pi}{\left(1-t^{\sum_{j=1}^k S_j D_j}\right)\left(1-t^{\sum_{j=1}^k S_j N_j}\right)}\left(\left[\text{gcd}\left(\sum_{j=1}^k S_j D_j, \left\{(S_1-i)D_1 + \sum_{l=2}^k S_lD_l\right\}_{i=1}^{S_1}, \{iN_1\}_{i=1}^{S_1}\right)\right] \right)\\
=&\dfrac{\Pi}{\left(1-t^{\sum_{j=1}^k S_j D_j}\right)\left(1-t^{\sum_{j=1}^k S_j N_j}\right)}(1-t)
\end{split}
\end{equation*}
where again we have assumed $S_1\geq 2$.
The Alexander polynomial is thus given by
$$\Delta_L(t) =(-1)^{-1+\sum_{j=1}^k S_j }t^{\frac\sigma 2} \dfrac{(1-t)\displaystyle\prod_{j=1}^k \left(1- t^{S_jN_jD_j +N_j \sum_{l=j+1}^k S_lD_l + D_j\sum_{k=1}^{j-1}S_lN_l}\right)^{S_j}}{\left(1-t^{\sum_{j=1}^k S_j D_j}\right)\left(1-t^{\sum_{j=1}^k S_j N_j}\right)}$$
where $$\sigma =-1+\sum_{j=1}^k S_j (N_j+D_j) - \sum_{j=1}^k S_j\left(S_jN_jD_j +N_j \sum_{l=j+1}^k S_lD_l + D_j\sum_{k=1}^{j-1}S_lN_l\right)$$

\section{Jones polynomials}
Jones, using representation theory and properties of Temperley-Lieb algebra, computed the Jones polynomial of the torus knot $K_{p, q}$ to be
$$V_{K_{p,q}}(t) = t^{\frac{(p-1)(q-1)}2} \dfrac{1-t^{p+1} - t^{q+1} +t^{p+q+1}}{1-t^2}$$
for coprime integers $p$ and $q$  \cite{Jones}. We now extend this to tst links. First we note that if $K^*$ denotes the mirror image of the knot $K$, $V_{K^*}(t) = V_K(t^{-1})$, hence
$$V_{K_{-p,q}}(t) = V_{K_{p,q}}(t^{-1}) = t^{-\frac{(p-1)(q-1)}2} \dfrac{1-t^{-p-1} - t^{-q-1} +t^{-p-q-1}}{1-t^{-2}}$$
since $K_{-p,q}$ is the mirror image of $K_{p,q}$ \cite{Murasugi}.

\subsection{Tst links $[\Phi^\mu (n(\tau), d(\tau), M)]$ with even $M$}
We consider the case where $M$ is even. Note that the link is the same as the torus link $K_{-\mu PN(\tau),\mu PD(\tau)}$ . We can thus use the technique of topological quantum field theory. First recall that the unnormalized Jones polynomial $\langle L \rangle$ of a given link $L$ is related to the vacuum expectation value of the Wilson lines in the $(2+1)$-dimensional Chern-Simons theory as follows \cite{Labastida}. First the partition function over the $3$-manifold $M$ is given by$$Z(M) = \int [DA]  \exp\left (\dfrac {ik}{8\pi\hbar} \int_M\text{Tr}  \left (A\wedge dA + \frac 23A\wedge A\wedge A\right)\right) $$
where $S_{CS} = \dfrac k{8\pi}\displaystyle \int_M\text{Tr}  \left (A\wedge dA + \frac 23A\wedge A\wedge A\right)$ is called the level-$k$ Chern-Simons action. The $\mu(L)$ components of $L$ are taken as Wilson lines and the partition function of $M$ with $L$ inside it is given by
$$Z(M; L) = \int [DA]  \exp \left(\dfrac i\hbar S_{CS}\right)\prod_{i=1}^{\mu(L)} \text{Tr}_R\mathcal{P}\exp \int_{L_i} A$$
where $L=\{L_1, L_2, ..., L_{\mu(L)}\}$ contains the components $L_i$, $W_R(C)=\text{Tr}_R\mathcal{P}\exp \displaystyle\int_C A$ is called the Wilson operator, $R$ is the $2$-dimensional representation of $SU(2)$ in the case of the usual Jones polynomial, and $\mathcal{P}$ is the path-ordering operator.

Then we can split the manifold $M$ containing $L$ into two parts $M_1$ and $M_2$ along a Riemann surface $\Sigma$ such that $M=M_1\cup_\Sigma M_2$ and no parts of $L$ appear in $\Sigma$. Suppose $L_1, L_2, ..., L_n$ lie in $M_1$ and the remaining components lie in $M_2$. Then we can define coordinates $z$ and $\bar{z}$ on the surface $\Sigma$, and the gauge connections $A_z$ and $A_{\bar{z}}$. $A_z$ and $A_{\bar{z}}$ satisfy a conjugate commutator relation. Choosing a metric $\gamma$ on $\Sigma$, we also define the wave functional associated with $M_1$ to be
$$\zeta_1[A_{\bar{z}}] = \int[DA]\prod_{i=1}^n W_R(L_i)\exp \left( \dfrac{ik}\hbar S_{CS} -\dfrac k{2\pi \hbar} \int_\Sigma \frac i2 dzd\bar{z} \sqrt\gamma \gamma^{z\bar{z}} \text{Tr}(A_zA_{\bar{z}})\right) $$
The functional corresponding to $M_2$ is given by
$$\zeta_2[A_{\bar{z}}] = \int[DA]\prod_{i=n+1}^{\mu(L)} W_R(L_i)\exp \left( \dfrac{ik}\hbar S_{CS} -\dfrac k{2\pi \hbar} \int_\Sigma \frac i2 dzd\bar{z} \sqrt\gamma \gamma^{z\bar{z}} \text{Tr}(A_zA_{\bar{z}})\right) $$
Then the inner product defined by
$$\langle \zeta_1|\zeta_2 \rangle = \int[DA_zDA_{\bar{z}}]\exp \left( \dfrac k{\pi \hbar} \int_\Sigma \frac i2 dzd\bar{z} \sqrt\gamma \gamma^{z\bar{z}} \text{Tr}(A_zA_{\bar{z}})\right)\overline{\zeta_1[A_{\bar{z}}]}\zeta_2[A_{\bar{z}}]  $$
gives both the correct commutator relation and the partition function, i.e.
$$Z(M;L) = \langle \zeta_1|\zeta_2 \rangle$$
Following \cite{Witten}, we consider the manifolds $M_1, M_2$ without the link components. Let the corresponding wave functionals be $\zeta_1'$ and $\zeta_2'$ respectively, then
\begin{equation*}
\begin{split}
Z(M) &= \langle \zeta_1'|\zeta_2' \rangle \\
Z(M_1; L_1, ..., L_n) &= \langle \zeta_1|\zeta_2' \rangle\\
Z(M_2; L_{n+1}, ..., L_{\mu(L)}) &= \langle \zeta_1'|\zeta_2 \rangle
\end{split}
\end{equation*}
We further note that
\begin{equation*}
\begin{split}
Z(M;L) Z(M) =& \langle \zeta_1|\zeta_2 \rangle\langle \zeta_1'|\zeta_2' \rangle \\
                       =& \int[DA_zDA_{\bar{z}}]\exp \left( \dfrac k{\pi \hbar} \int_\Sigma \frac i2 dzd\bar{z} \sqrt\gamma \gamma^{z\bar{z}} \text{Tr}(A_zA_{\bar{z}})\right)\overline{\zeta_1[A_{\bar{z}}]}\zeta_2[A_{\bar{z}}]\\
                       &\times \int[DA_zDA_{\bar{z}}]\exp \left( \dfrac k{\pi \hbar} \int_\Sigma \frac i2 dzd\bar{z} \sqrt\gamma \gamma^{z\bar{z}} \text{Tr}(A_zA_{\bar{z}})\right)\overline{\zeta_1'[A_{\bar{z}}]}\zeta_2'[A_{\bar{z}}]\\
                       =& \int\int[DA_zDA_{\bar{z}}]\exp \left( \dfrac k{\pi \hbar} \int_\Sigma i dzd\bar{z} \sqrt\gamma \gamma^{z\bar{z}} \text{Tr}(A_zA_{\bar{z}})\right)\overline{\zeta_1[A_{\bar{z}}]}\zeta_2[A_{\bar{z}}]\overline{\zeta_1'[A_{\bar{z}}]}\zeta_2'[A_{\bar{z}}]\\
                        =& \int\int[DA_zDA_{\bar{z}}]G^2\overline{\zeta_1'[A_{\bar{z}}]}\zeta_2[A_{\bar{z}}]\overline{\zeta_1[A_{\bar{z}}]}\zeta_2'[A_{\bar{z}}]\\
                         =& \int [DA_zDA_{\bar{z}}]G\overline{\zeta_1'[A_{\bar{z}}]}\zeta_2[A_{\bar{z}}]\int [DA_zDA_{\bar{z}}]G\overline{\zeta_1[A_{\bar{z}}]}\zeta_2'[A_{\bar{z}}]\\
                         =&\langle \zeta_1'|\zeta_2 \rangle \langle \zeta_1|\zeta_2'\rangle \\
                         =&Z(M_1; L_1, ..., L_n) Z(M_2; L_{n+1}, ..., L_{\mu(L)})
\end{split}
\end{equation*}
where $G=\exp \left( \dfrac k{\pi \hbar} \displaystyle\int_\Sigma \dfrac i2 dzd\bar{z} \sqrt\gamma \gamma^{z\bar{z}} \text{Tr}(A_zA_{\bar{z}})\right)$. 
Rearranging gives
$$
\dfrac{Z(M;L)}{Z(M)} = \dfrac{Z(M_1; L_1, ..., L_n) }{Z(M)}\dfrac{Z(M_2; L_{n+1}, ..., L_{\mu(L)})}{Z(M)}
$$

Hence we are motivated to define the unnormalized Jones polynomial $\langle L \rangle$ of the link $L$ to be
$$\langle L \rangle = \dfrac{Z(M;L)}{Z(M)}$$
so that 
\begin{equation}
\label{multiply}
\langle L\rangle=\langle  L_1, ..., L_n \rangle \langle L_{n+1}, ..., L_{\mu(L)}\rangle
\end{equation}

Now the normalized Jones polynomial is given by
$$V_L(t) = \dfrac{\langle L \rangle}{\langle O \rangle}$$
where $O$ denotes the unknot. From \cite{Witten}, we also have
$$\langle O \rangle = \dfrac{t-t^{-1}}{t^{\frac 12} - t^{-\frac 12}} = t^{-\frac 12} \dfrac{1-t^2}{1-t} = t^{\frac 12} \dfrac{1-t^{-2}}{1-t^{-1}}$$

To consider the case where $L =[\Phi^\mu (n(\tau), d(\tau), M)]$ with even $M$, we first compute $\langle K_{-N, D}\rangle$:
\begin{equation*}
\begin{split}
\langle K_{-N, D}\rangle &= V_{K_{-N, D}} (t) t^{-\frac 12} \dfrac{1-t^2}{1-t}\\
&=  t^{-\frac{(N-1)(D-1)}2} \dfrac{1-t^{-N-1} - t^{-D-1} +t^{-N-D-1}}{1-t^{-2}} t^{\frac 12} \dfrac{1-t^{-2}}{1-t^{-1}} \\
&=t^{\frac{1-(N-1)(D-1)}2} \dfrac{1-t^{-N-1} - t^{-D-1} +t^{-N-D-1}}{1-t^{-1}}\\
&=t^{\frac{N+D-ND}2} \dfrac{1-t^{-N-1} - t^{-D-1} +t^{-N-D-1}}{1-t^{-1}}
\end{split}
\end{equation*}
We then consider $(S-1) $ solid tori $U_1,U_2,\cdots, U_{S-1} $ that share a common longitudinal circle such that $U_{S-1} \subset \cdots \subset U_2 \subset U_1$. We then consider another set of $S$ solid tori $V_1,V_2,\cdots, V_S $ that share the same  longitudinal circle such that $\partial V_i \subset U_{i-1} \backslash U_i $ for $i\geq 2$ and $\partial V_1 \subset S^3\backslash U_1$ 	and place each component $K_{-N,D} $ on $\partial V_i$. Then we can split $M=S^3$ into the spaces $M_1 = \overline {S^3\backslash U_1} $ and $M_2 = \overline{U_1} $ so that $\Sigma = \partial U_1$ and $M=M_1\cup_\Sigma M_2$, and apply equation  (\ref{multiply}) to give
$$\langle L\rangle = t^{\frac{N+D-ND}2} \dfrac{1-t^{-N-1} - t^{-D-1} +t^{-N-D-1}}{1-t^{-1}} \langle L'\rangle $$
where $L'$ contains the remaining $(S-1)$ components.
We then further split $\overline{U_1} $ into $\overline{U_2} $ and $\overline{U_1\backslash U_2} $ and apply (\ref{multiply}) again, until we are left with the last component.
Thus we have
\begin{equation*}
\begin{split}
\langle K_{-SN,SD} \rangle &= \left(t^{\frac{N+D-ND}2} \dfrac{1-t^{-N-1} - t^{-D-1} +t^{-N-D-1}}{1-t^{-1}} \right) ^S\\
&=t^{\frac{S(N+D-ND)}2}\left(\dfrac{1-t^{-N-1} - t^{-D-1} +t^{-N-D-1}}{1-t^{-1}} \right) ^S
\end{split}
\end{equation*} 
The normalized Jones polynomial is given by
\begin{equation*}
\begin{split}
V_{[\Phi^\mu (n(\tau),d(\tau),M) ]} (t) &= \dfrac{\langle K_{-SN,SD} \rangle} {\langle O\rangle}\\
&=\dfrac{t^{\frac{S(N+D-ND)}2}\left(\dfrac{1-t^{-N-1} - t^{-D-1} +t^{-N-D-1}}{1-t^{-1}} \right) ^S}{t^{\frac 12} \dfrac{1-t^{-2}}{1-t^{-1}}}\\
&=t^{\frac{S(N+D-ND)-1}2}\dfrac{(1-t^{-N-1} - t^{-D-1} +t^{-N-D-1})^S}{(1-t^{-1})^S(1+t^{-1})}
\end{split}
\end{equation*} 

\subsection{Tst links $[\Phi^\mu (n(\tau), d(\tau), M)]$ with odd $M$}
We now turn to the case where $M$ is odd. We can simply set $M_2 = \overline{U}$ to be a solid torus containing the trivial knot and $M_1=\overline {S^3 \backslash U} $ to contain the torus link. Then $\Sigma = \partial U$ and we can apply equation (\ref{multiply}) to give
$$\langle [\Phi^\mu (n(\tau), d(\tau), M)] \rangle= \langle K_{-SN,SD}\rangle \langle O \rangle$$
The normalized Jones polynomial is given by
$$V_{[\Phi^\mu (n(\tau), d(\tau), M)]} (t) = \dfrac{\langle K_{-SN,SD}\rangle \langle O \rangle} {\langle O \rangle}=t^{\frac{S(N+D-ND)}2}\left(\dfrac{1-t^{-N-1} - t^{-D-1} +t^{-N-D-1}}{1-t^{-1}} \right) ^S$$
\subsection{Stst links $[\Phi_{j=1}^k(n(\tau_j), d(\tau_j), M_j,\mu_j)]$ with even $M_k$}
We now turn to the stst link $[\Phi_{j=1}^k(n(\tau_j), d(\tau_j), M_j,\mu_j)]$. We write $N(\tau_j)=N_j $ , $D(\tau_j)=D_j$ and $\mu_jP_j = S_j$. For even $M_k$, we consider $k$ solid tori with the same longitudinal circle, $V_k\subset \cdots \subset V_2\subset V_1$ and another $(k-1)$ solid tori with also the same longitudinal circle, $U_{k-1}\subset \cdots \subset U_2\subset U_1$  and place $\partial V_1$ inside $\overline{S^3\backslash U_1}$ and $\partial V_i$ inside $\overline{U_{i-1}\backslash U_i}$ for $i\geq 2$. On each $\partial V_j$, we place the torus link $[\Phi^{\mu_j}(n(\tau_j), d(\tau_j),M_j)]$. We now repeatedly apply equation(\ref{multiply}) so that
\begin{equation*}
\begin{split}
V_{[\Phi_{j=1}^k(n(\tau_j), d(\tau_j), M_j,\mu_j)]}(t) &= \dfrac{\displaystyle\prod_{j=1}^k t^{\frac{S_j(N_j+D_j-N_jD_j)}2}\left(\dfrac{1-t^{-N_j-1} - t^{-D_j-1} +t^{-N_j-D_j-1}}{1-t^{-1}} \right) ^{S_j}}{t^{\frac 12} \dfrac{1-t^{-2}}{1-t^{-1}}}\\
&= t^{\frac {-1+\sum_{j=1}^k S_j(N_j+D_j-N_jD_j)}2}\dfrac{\displaystyle\prod_{j=1}^k (1-t^{-N_j-1} - t^{-D_j-1} +t^{-N_j-D_j-1})^{S_j}}{(1-t^{-1})^{\sum_{j=1}^kS_j}(1+t^{-1})}
\end{split}
\end{equation*}
\subsection{Stst links $[\Phi_{j=1}^k(n(\tau_j), d(\tau_j), M_j,\mu_j)]$ with odd $M_k$}
We now turn to the case where $M_k$ is odd. This is similar to the case where $M_k$ is even, but we introduce one more solid torus $U_k$ that shares the same longitudinal circle with other $U$'s such that the trivial knot is placed inside $U_k$. We can now apply equation (\ref{multiply}) so that
$$V_{[\Phi_{j=1}^k(n(\tau_j), d(\tau_j), M_j,\mu_j)]}(t) = t^{\frac {\sum_{j=1}^k S_j(N_j+D_j-N_jD_j)}2}\dfrac{\displaystyle\prod_{j=1}^k (1-t^{-N_j-1} - t^{-D_j-1} +t^{-N_j-D_j-1})^{S_j}}{(1-t^{-1})^{\sum_{j=1}^kS_j}} $$

\section*{Conclusions}
In the paper we introduced stst links, computed their braid words, Alexander and Jones polynomials. Further research directions include finding other knot or link invariants such as unknotting numbers, hyperbolic volumes, HOMFLY polynomials, Kontsevich integrals and other Vassiliev invariants. Further generalizations also include satellite links built from stst links, and virtual knots. These could be investigated in the future.


\end{document}